\documentclass{article}[12pt]

\usepackage{latexsym}
\usepackage{amsmath}
\usepackage{amssymb}

\begin{document}

\title{EXPLICITLY SOLVABLE SYSTEMS OF FIRST-ORDER DIFFERENCE EQUATIONS WITH
\ HOMOGENEOUS POLYNOMIAL RIGHT-HAND SIDES}

\author{Francesco Calogero$^{a,b}$\thanks{e-mail: francesco.calogero@roma1.infn.it}
\thanks{e-mail: francesco.calogero@uniroma1.it}
 , Farrin Payandeh$^c$\thanks{e-mail: farrinpayandeh@yahoo.com}
 \thanks{e-mail: f$\_$payandeh@pnu.ac.ir}}

\maketitle   \centerline{\it $^{a}$Physics Department, University of
Rome "La Sapienza", Rome, Italy}

\maketitle   \centerline{\it $^{b}$INFN, Sezione di Roma 1}

\maketitle

\maketitle   \centerline{\it $^{c}$Department of Physics, Payame
Noor University, PO BOX 19395-3697 Tehran, Iran}

\maketitle

In this short paper we identify \textit{special} systems of (an \textit{%
arbitrary} number) $N$ of first-order Difference Equations with \textit{%
nonlinear homogeneous} polynomials of \textit{arbitrary} degree $M$ in their
right-hand sides, which feature \textit{very simple explicit }solutions. A
novelty of these findings is to consider \textit{special} systems
characterized by \textit{constraints} involving both their parameters and
their initial data.

\bigskip

The general system of an \textit{arbitrary} number $N$ of first-order
Difference Equations (DEs) with \textit{homogeneous} \textit{polynomials} of
\textit{arbitrary} degree $M$ on their right-hand sides reads as follows:
\begin{eqnarray}
\tilde{z}_{n}\left( s\right) &=&\sum_{m_{\ell }}{}^{\left( M\right) }\left\{
c_{nm_{1}m_{2}\cdot \cdot \cdot m_{N}}\left[ z_{1}\left( s\right) \right]
^{m_{1}}\left[ z_{2}\left( s\right) \right] ^{m_{2}}\cdot \cdot \cdot \left[
z_{N}\left( s\right) \right] ^{m_{N}}\right\} ~,  \notag \\
n &=&1,2,...,N~.  \label{1}
\end{eqnarray}

\textbf{Notation. }Above and hereafter $s$ is the \textit{discrete}
independent variable taking all nonnegative integer values, $s=0,1,2,...$;
the $N$ dependent variables $z_{n}\left( s\right) $ with $n=1,2,...,N$ are
(possibly \textit{complex}) numbers, and ascertaining their $s$-evolution
from the set of $N$ \textit{initial data} $z_{n}\left( 0\right) $ is our
main task; the symbol $\tilde{z}_{n}\left( s\right) $ denotes{}the
forward-shifted dependent variable,
\begin{equation}
\tilde{z}_{n}\left( s\right) \equiv z_{n}\left( s+1\right) ~,~~~n=1,2,...,N~;
\label{ztilde}
\end{equation}%
the symbol $\sum_{m_{\ell }}{}^{\left( M\right) }$ denotes the sum running
over \textit{all nonnegative} values of the $N$ \textit{nonnegative integer}
parameters (indices and exponents) $m_{\ell }$ subject to the restrictions
\begin{equation}
m_{\ell }\geq 0~,~~~\sum_{\ell =1}^{N}\left( m_{\ell }\right) =M~,  \label{M}
\end{equation}%
implying that the polynomials in $N$ variables $z_{n}\left( s\right) $ in
the right-hand sides of the $N$ DEs (\ref{1}) are \textit{all homogeneous}
of degree $M$, being characterized by the $s$-independent coefficients $%
c_{nm_{1}m_{2}\cdot \cdot \cdot m_{N}}$. $\blacksquare $

The findings reported in this paper are the extension to \textit{discrete}
time of the somewhat analogous results for systems of first-order Ordinary
Differential Equations (ODEs) reported in \cite{CP2021}; indeed, its
presentation occasionally follows \textit{verbatim} the text of \cite{CP2021}%
. Simple as they are, they are to the best of our knowledge \textit{new},
being based on a somewhat \textit{unconventional} approach: to identify
\textit{explicitly solvable} cases of the system (\ref{1}) by introducing
\textit{constraints} involving, in addition to the coefficients $%
c_{nm_{1}m_{2}\cdot \cdot \cdot m_{M}}$, also the initial data $z_{n}\left(
0\right) $ (which, in applicative contexts, may play the role of \textit{%
control} elements, allowing to manipulate the time evolution of the system).

\textbf{Remark 1}. In this paper we focus on systems with \textit{%
homogeneous polynomial} right-hand sides, see (\ref{1}); but clearly---as
pointed out by the extension \cite{PZ2021} of the results of \cite{CP2021}%
---these findings can be extended to \textit{more general} \textit{%
homogeneous} functions than polynomials (see \cite{CP2021bis} \cite%
{CP2021ter}). $\blacksquare $

Our main result is the following

\textbf{Proposition}. The system (\ref{1}) features the special solution
\begin{subequations}
\label{Solznt}
\begin{equation}
z_{n}\left( s\right) =z_{n}\left( 0\right) \left[ z_{N}\left( 0\right) %
\right] ^{M^{s}-1}Z^{\left( M^{s}-1\right) /\left( M-1\right)
}~,~~~n=1,2,...,N~,  \label{znt}
\end{equation}%
provided there hold the following $N$ explicit \textit{algebraic constraints}
on the \textit{a priori arbitrary }parameter $Z$, the coefficients $%
c_{nm_{1}m_{2}\cdot \cdot \cdot m_{M}}$ and the $N$ initial data $%
z_{n}\left( 0\right) $:
\begin{equation}
Z=\left( r_{n}\right) ^{-1}\sum_{m_{\ell }}{}^{\left( M\right) }\left\{
c_{nm_{1}m_{2}\cdot \cdot \cdot m_{N}}\prod\limits_{\ell =1}^{N-1}\left[
\left( r_{\ell }\right) ^{m_{\ell }}\right] \right\} ~,~~~n=1,2,...,N~,
\label{Kzn0}
\end{equation}%
where (here and hereafter)%
\begin{equation}
r_{n}\equiv z_{n}\left( 0\right) /z_{N}\left( 0\right) ~.~\blacksquare
\label{rn}
\end{equation}

\textbf{Remark 2}. The proof that (\ref{Solznt}) satisfies the system of
ODEs (\ref{1}) is elementary: just insert (\ref{znt}) in (\ref{1}) and
verify that, thanks to (\ref{M}) and (\ref{Kzn0}), the $N$ DEs (\ref{1}) are
satisfied. $\blacksquare $

\textbf{Remark 3}. Note than only the ratios of the $N$ initial data play a
role in the constraints (\ref{Kzn0}).

\textbf{Remark 4}. The system of $N$ \textit{algebraic} equations (\ref{Kzn0}%
) generally determines---for any given assignment of the \textit{a priori
arbitrary} coefficients $c_{nm_{1}m_{2}\cdot \cdot \cdot m_{N}}$---$N$ out
of the $N+1$ quantities $Z$ and $z_{n}\left( 0\right) $ (of the latter, only
their ratio); but it is also possible to select \textit{ad libitum} $N$
elements out of the \textit{complete} set of data $Z$, $c_{nm_{1}m_{2}\cdot
\cdot \cdot m_{M}}$ and $z_{n}\left( 0\right) $ (of the latter, only their
ratio), and to then consider these \textit{selected} elements as those to be
determined---by the $N$ conditions (\ref{Kzn0})---in terms of the remaining
\textit{arbitrarily assigned} elements in the \textit{complete} set of these
data. If one chooses to satisfy these $N$ conditions by solving the $N$
equations (\ref{Kzn0}) for $N$ of the coefficients $c_{nm_{1}m_{2}\cdot
\cdot \cdot m_{M}}$---or for the parameter $Z$ and $N-1$ of the coefficients
$c_{nm_{1}m_{2}\cdot \cdot \cdot m_{M}}$---then this task can be generally
performed \textit{explicitly}, since the relevant \textit{algebraic}
equations to be solved are then \textit{linear} in the unknown quantities;
otherwise these determinations require the solution of \textit{nonlinear }%
equations, a task which can be performed \textit{explicitly} only rarely in
an \textit{algebraic} setting; but which can generally be performed, with
\textit{arbitrary} approximation, in a \textit{numerical} context. $%
\blacksquare $

\textbf{Example}. Assume for instance $N=2$ and $M=4$, so that the system (%
\ref{1}) reads as follows (note below the notational simplification):
\end{subequations}
\begin{equation}
\tilde{z}_{n}\left( s\right) =\sum_{m=0}^{4}c_{nm}\left[ z_{1}\left(
s\right) \right] ^{4-m}\left[ z_{2}\left( s\right) \right] ^{m}~,~~~n=1,2~,
\end{equation}%
featuring $2$ dependent variables $z_{n}\left( s\right) $ and $10$ \textit{a
priori arbitrary} coefficients $c_{nm}$ ($n=1,2$; $m=0,1,2,3,4$). Then the
solution (\ref{znt}) reads as follows:
\begin{subequations}
\label{Solznt}
\begin{equation}
z_{n}\left( s\right) =z_{n}\left( 0\right) \left[ z_{2}\left( 0\right) %
\right] ^{4^{s}-1}Z^{\left( 4^{s}-1\right) /3}~,~~~n=1,2,...,N~,
\end{equation}%
and the $2$ conditions (\ref{Kzn0}) read as follows:%
\begin{equation}
Z=r^{-1}\sum_{m=0}^{2}\left( c_{1m}r^{m}\right) =\sum_{m=0}^{2}\left(
c_{2m}r^{m}\right) ~,
\end{equation}%
with $r\equiv z_{1}\left( 0\right) /z_{2}\left( 0\right) $, namely%
\begin{equation}
\sum_{m=0}^{2}\left[ \left( c_{1m}-c_{2m}r\right) r^{m}\right] =0~.
\end{equation}

These $2$ algebraic constraints can of course be \textit{explicitly} solved
for any $2$ of the $10$ coefficients $c_{nm}$ in terms of the other $8$
coefficients $c_{nm}$ and of the $2$ \textit{arbitrary} data $Z$ and the
ratio $r\equiv z_{1}\left( 0\right) /z_{2}\left( 0\right) $; or
alternatively for $Z$ and only $1$ of the $10$ coefficients $c_{nm}$ in
terms of the other $9$ coefficients $c_{nm}$ and of the ratio $r\equiv
z_{1}\left( 0\right) /z_{2}\left( 0\right) $; with many other possibilities
left to the imagination of the interested reader. $\blacksquare $

\textbf{Final Remark}. As already noted above, the mathematics behind the
results reported above is rather \textit{elementary}. Yet these findings do
not seem to have been advertised so far, while their \textit{applicable}
potential is clearly vast; so---especially among \textit{applied}
mathematicians and \textit{practitioners} of the various scientific
disciplines where systems such as those discussed above play a role---a
wider knowledge of them seems desirable; for instance via their inclusion in
standard compilations of \textit{solvable} equations such as those collected
in the website EqWorld. $\blacksquare $

\textbf{Acknowledgements}. We like to acknowledge with thanks $2$ grants,
facilitating our collaboration---mainly developed via e-mail exchanges---by
making it possible for FP to visit twice the Department of Physics of the
University of Rome "La Sapienza": one granted by that University, and one
granted jointly by the Istituto Nazionale di Alta Matematica (INdAM) of that
University and by the International Institute of Theoretical Physics (ICTP)
in Trieste in the framework of the ICTP-INdAM "Research in Pairs" Programme.%
\textbf{\ }Finally, we also\ like to thank Fernanda Lupinacci who, in these
difficult times---with extreme efficiency and kindness---facilitated all the
arrangements necessary for the presence of FP with her family in Rome.

\end{subequations}

\end{document}